\catcode`\@=11
\magnification 1200
\hsize=140mm \vsize=200mm
\hoffset=-4mm \voffset=-1mm
\pretolerance=500 \tolerance=1000 \brokenpenalty=5000

\catcode`\;=\active
\def;{\relax\ifhmode\ifdim\lastskip>\z@
\unskip\fi\kern.2em\fi\string;}

\catcode`\:=\active
\def:{\relax\ifhmode\ifdim\lastskip>\z@\unskip\fi
\penalty\@M\ \fi\string:}

\catcode`\!=\active
\def!{\relax\ifhmode\ifdim\lastskip>\z@
\unskip\fi\kern.2em\fi\string!}\catcode`\?=\active
\def?{\relax\ifhmode\ifdim\lastskip>\z@
\unskip\fi\kern.2em\fi\string?}

\frenchspacing

\newif\ifpagetitre        \pagetitretrue
\newtoks\hautpagetitre    \hautpagetitre={\hfil}
\newtoks\baspagetitre     \baspagetitre={\hfil}
\newtoks\auteurcourant    \auteurcourant={\hfil}
\newtoks\titrecourant     \titrecourant={\hfil}
\newtoks\hautpagegauche   \newtoks\hautpagedroite
\hautpagegauche={\hfil\tensl\the\auteurcourant\hfil}
\hautpagedroite={\hfil\tensl\the\titrecourant\hfil}

\newtoks\baspagegauche
\baspagegauche={\hfil\tenrm\folio\hfil}
\newtoks\baspagedroite
\baspagedroite={\hfil\tenrm\folio\hfil}

\headline={\ifpagetitre\the\hautpagetitre
\else\ifodd\pageno\the\hautpagedroite
\else\the\hautpagegauche\fi\fi}

\footline={\ifpagetitre\the\baspagetitre
\global\pagetitrefalse
\else\ifodd\pageno\the\baspagedroite
\else\the\baspagegauche\fi\fi}

\font\twbf=cmbx12\font\sc=cmcsc10

\def\date{le\ {\the\day}\ 
\ifcase\month\or Janvier\or \F\'evrier\or Mars\or Avril
\or Mai\or Juin\or Juillet\or Ao\^ut\or Septembre
\or Octobre\or Novembre\or D\'ecembre\fi\ {\the\year}}

\def\cf{{\it cf.\/}\ }    
 \def\up#1{\raise 1ex\hbox{\sevenrm#1}}
\def\cqfd{\unskip\kern 6pt\penalty 500
\raise -2pt\hbox{\vrule\vbox to 10pt{\hrule width 4pt\vfill
\hrule}\vrule}\par}
\catcode`\@=12
\def \bg {\bigskip \goodbreak}

\def\ref#1&#2&#3&#4&#5\par{\par{\leftskip = 5em {\noindent
\kern-5em\vbox{\hrule height0pt depth0pt width
5em\hbox{\bf[\kern2pt#1\unskip\kern2pt]\enspace}}\kern0pt}
{\sc\ignorespaces#2\unskip},\
{\rm\ignorespaces#3\unskip}\
{\sl\ignorespaces#4\unskip\/}\
{\rm\ignorespaces#5\unskip}\par}}

\def\exo#1{\goodbreak\vskip 12pt plus 20pt minus 2pt
\line{\noindent\hss\bf
\uppercase\expandafter{\romannumeral#1}\hss}\nobreak\vskip
12pt }
\def \titre#1\par{\null\vskip
1cm\line{\hss\vbox{\twbf\halign
{\hfil##\hfil\crcr#1\crcr}}\hss}\vskip 1cm}

\def \frac#1#2{{{#1}\over {#2} }}

\def \comp{ \;{}^{ {}_\vert }\!\!\!{\rm C}   }

\def \nat{ { {\rm I}\!{\rm N}} }

\def \rat{ {\rm Q}\kern-.65em {}^{{}_/ }}

\def \adh {\overline}

\def\N#1{\left\Vert #1\right\Vert }

\def\dess#1by#2(#3){\vbox to #2{
\hrule width #1 height 0pt depth 0pt
\vfill\special{picture #3}
}}

\baselineskip=18pt
\def\adh{\overline}
Analyse Fonctionnelle

\centerline {Espace de Hilbert d'op\'erateurs}
\centerline {et}
\centerline {Interpolation complexe}
\centerline {par}
\centerline {Gilles Pisier}
Note pr\'esent\'ee par Alain Connes
\bg 
\bg

\bg 
\bg 

\underbar {R\'esum\'e}. Soit $H$ un espace hilbertien de dimension infinie. On montre qu'il existe un sous-espace ferm\'e de $B(H)$ qui
est hilbertien
 et compl\`etement isom\'etrique \`a son antidual au sens de la th\'eorie des espaces
d'op\'erateurs d\'evelop\'ee r\'ecemment par Blecher-Paulsen et Effros-Ruan. De plus cet
espace est unique \`a une isom\'etrie compl\`ete pr\`es. On le note $OH$. Cet espace a de
nombreuses propri\'et\'es remarquables, en particulier pour l'interpolation
complexe. \bg  \bg 
\underbar {Abstract}. Let $H$ be an infinite dimensional Hilbert space. We show
that there exists a subspace of $B(H)$ which is isometric
to $\ell_2$ and completely isometric to its antidual in the
sense of the theory of operator spaces recently developed
by Blecher-Paulsen and Effros-Ruan. Moreover this space is
unique up to a complete isometry. We denote it by $OH$.
This space has several remarkable properties in particular
with respect to the complex interpolation method. \bg  \bg 

\underbar {Abridged English Version}. 

In the theory of operator spaces as recently developed in [2,6] several
examples of Hilbertian operator spaces play an important role for instance the
spaces $C = B(\comp,\ell_2)$ and $R = B(\ell_2,\comp)$. However in general a
Hilbertian operator space $E$ needs not be completely isomorphic to its
``operator dual" $E^*$ as defined in [2,6]. Let $H$ the Hilbert space $\ell_2$.
One of the main results of this announcement is the observation that there is a
subspace $OH \subset B(H)$ which is isometric to $\ell_2$ and such that the
natural linear isomorphism from $OH$ to $\adh {OH^*}$ is a complete isometry.
Moeover this space $OH$ is unique up to a complete isometry.

Let $E,F$ be two closed subspaces of $B(H)$. We denote by $E \otimes_{\min} F$
their minimal (or spatial) tensor product in $B(H \otimes H)$. If $H,K$ are
Hilbert spaces we denote by $H \otimes K$ their Hilbertian tensor product.

Let $(T_n)_{n \geq 0}$ be an orthomormal basis of $OH$ and let $x_n$
 be a finite sequence in $B(H)$. We have then
$$\N {\sum T_n \otimes x_n}_{OH \otimes_{\min} B(H)} = \N {\sum x_n \otimes
\adh{x}_n}_{B(H \otimes \overline {H})}^{1/2}.$$
The proof of the preceding results makes crucial use of an inequality due to U.
Haagerup ([7], lemma 2.4) which seems to be the appropriate operator version of
the Cauchy-Schwarz inequality.

In the second part, we study the complex interpolation method between couples
of operator spaces. Let $(E_0,E_1)$ be compatible couple of Banach spaces in
the sense of interpolation theory (cf. [1]). Assume that $E_0$ and $E_1$ are
both equipped with an operator space structure in the sense of [2,6]. Then the
complex interpolation spaces $E_\theta = (E_0,E_1)_\theta$ and $E^\theta =
(E_0,E_1)^\theta$ in the sense of [1] can be equipped with an operator space
structure by defining
$$\forall~n \geq 1~~~~~~~M_n(E_\theta) = (M_n(E_0),M_n(E_1))_\theta~{\rm
and}~M_n(E^\theta) = (M_n(E_0),M_n(E_1))^\theta.$$
The proof that these are indeed operator spaces uses Ruan's criterion [11].

Let $E_0 \cap E_1$ be the intersection equipped with the norm $\N {x} = \max
\{\N {x}_{E_0},\N {x}_{E_1}\}$ for all $x$ in  $E_0 \cap E_1$. Similarly for the
sum $E_0 + E_1$ we define $\N {x}_{E_0 + E_1} = \inf \{\N {y}_{E_0} + \N
{z}_{E_1} |x = y + z\}$.

 We denote by 
$E_0 \oplus_\infty E_1$ (resp. $E_0 \oplus_1 E_1$) the
direct sum equipped with the norm $\|(x,y)\|=
\max(\|x\|,\|y\|)$ (resp. $\|x\|+\|y\|)$. We equip these
spaces with the operator space structure defined by setting
$M_n(E_0 \oplus_\infty E_1)=M_n(E_0 )\oplus_\infty
M_n(E_1)$
(resp. $M_n(E_0 \oplus_1 E_1)$ is equipped with the norm
induced by the space 
$cb(E_0^*\oplus_{\infty} E_1^*,M_n)$).

 We can also equip the spaces $E_0 \cap E_1$
and $E_0 + E_1$ with a natural operator space structure. We
identify   $E_0 \cap E_1$ with the subspace
$\Delta=\{(x,x)\}$ in $E_0 \oplus_\infty E_1$ and identify
$E_0 + E_1$ with the    operator  space structure of the
quotient $(E_0 \oplus_1 E_1)/\{(x,-x)\}$.
 Since the column and row Hilbert spaces $C$
and $R$ are naturally both linearly identifiable with $\ell_2$, we may view
$(C,R)$ as a compatible interpolation couple and we obtain the operator spaces
$C \cap R, C + R$ and $(C,R)_\theta$. (The operator space structure of $C \cap
R$
is studied in [8].) We have completely bounded inclusions (actually complete
contractions) $$C \cap R \to OH \to C+R.$$ More generally let $E$ be any
operator space and let $v : E \to \ell_2$ be a bounded injective map with dense
range so that the mapping  $J = \adh {v^*} v : E \to \adh {E^*}$ is a continuous
injection allowing us to consider $(E,\adh {E^*})$ as a compatible couple. Then
the space $(E,\adh{E^*})_{1/2}$ is completely isometric to $OH$. In particular
$OH$ is completely isometric to $(R,C)_{1/2}$. We also study the Haagerup tensor
product (cf. [2,6]) between two operator spaces $E,F$. We denote it by $E
\otimes_h F$. The results of [9] imply that if $(E_0,E_1)$ and $(F_0,F_1)$ are
compatible couples and if $E_\theta = (E_0,E_1)_\theta,F_\theta =
(F_0,F_1)_\theta$ then $E_\theta \otimes_h F_\theta$ is completely isometric to
the interpolation space $(E_0 \otimes_h F_0,E_1 \otimes_h F_1)_\theta$.

Among various applications, we prove that a von Neumann algebra $M \subset
B(H)$ is injective if (and only if) there is a completely bounded linear
projection $P : B(H) \to M$. Actually, it suffices that $P$ satisfies for some
constant $C$ for all finite sequences $x_1,...,x_n$ in $B(H)$
$$\N {\sum P(x_i)^* P(x_i)} \leq C^2 \N {\sum x_i^* x_i}~{\rm and}$$
$$\N {\sum P(x_i) P(x_i)^*} \leq C^2 \N {\sum x_i x_i^*}.$$

 Let $E,F$ be operator spaces and let $u\in E\otimes F$. We define
 $$oh(u) = \inf\left\{\left\|\sum x_i\otimes \bar
 x_i\right\|^{1/2}_{E\otimes \overline E} \left\|\sum y_i\otimes \bar
 y_i\right\|^{1/2}_{F\otimes \overline F}\right\} $$
 where the infimum runs over all the factorizations of $u$ of the form
 $$u = \sum^n_{i=1} x_i \otimes y_i.$$
 It is very easy to check that this is a norm on $E\otimes F$.
 Equivalently, $oh(u)$ is the ``norm of factorization through $OH$'' of the
 linear operator $u_1\colon \ E^*\to F$ (or $u_2\colon  \ F^*\to E$)
 canonically associated to $u$. Indeed, by the preceding  results   we
have
 $$oh(u) = \inf\{\|A\|_{cb} \|B\|_{cb}\} $$
 where the infimum runs over all the possible factorizations $u_1 = AB$
 with $B\colon \ E^*\to OH$    and $A\colon \ OH\to F$ both completely bounded.
  We
 will denote by $E\otimes _{oh} F$ the completion of $E\otimes F$ for the
 above norm. This gives us a Banach space structure on
$E\otimes_{oh} F$.

Let $H,K$ be Hilbert spaces.
 Consider the linear mapping  $$\Phi\colon \ B(H)\otimes B(K)\to cb(B(H,K),
B(H,K))$$ defined on the linear tensor
product $B(H)\otimes B(K)$ by the property that it  
   maps $x\otimes y$ to the operator $\Phi(x\otimes
y)\colon \ B(H,K)\to B(H,K)$ defined by $$\Phi(x\otimes y)(z) =
xzy.$$ In the finite
 dimensional case, i.e. if say $\dim H = \dim K = n$ it is known that $\Phi$
 defines a (complete) isometric isomorphism between $M_n\otimes_h M_n$ and
 $cb(M_n,M_n)$. From this it is easy to deduce that the ``same'' map $\Phi$
 defines a (complete) isometric isomorphism from $M^{op}_n \otimes_h
 M^{op}_n$ onto $cb(\overline{M^*_n}, \overline{M^*_n})$ where we identify
 as sets (for the purpose of interpolation) the linear spaces $M_n$ and
 $\overline{M^*_n}$ by the usual canonical linear isomorphism $i\colon \
 M_n\to \overline{M^*_n}$ defined by $i(x)(y) = tr(xy^*)$. Then we have
 
 \proclaim Corollary.   With the preceding identification the map $\Phi$
 defines an isometric isomorphism from $M_n\otimes_{oh} M_n$ onto
 $(cb(M_n,M_n)$, $cb(\overline{M^*_n},\overline{M^*_n}))_{1\over 2}$.
 
A similar statement can be given in the setting of the spaces considered in [5]

and [3] of c.b. maps which are $(M',N')$ modular where $M\subset B(H)$ and 
$N\subset B(K)$ are injective von Neumann algebras.

 We will study the tensor product $E\otimes _{oh} F$ and the space of operators
factoring through $OH$ in a forthcoming note. \vfill\eject

Le but de cette note est d'annoncer les r\'esultats principaux de notre article
[10].
Nous utiliserons la th\'eorie des espaces d'op\'erateurs de [2,6]. Dans cette
th\'eorie plusieurs espaces hilbertiens  d'op\'erateurs jouent un r\^ole important,
par exemple les espaces  $C = B(\comp,\ell_2)$ (espace de Hilbert en colonne) et
$R = B(\ell_2,\comp)$ (espace de Hilbert en ligne), mais, contrairement aux
espaces de Hilbert usuels, aucun des espaces utilis\'es jusqu'ici ne semble \^etre
canoniquement isomorphe \`a son antidual au sens de [2,6]. L'un de nos 
principaux r\'esultats est l'observation qu'il existe un (et un seul \`a isom\'etrie
compl\`ete pr\'es) espace hilbertien d'op\'erateurs  qui est compl\`etement isom\'etrique
\`a son antidual.

 Soient $E,F$ deux espaces d'op\'erateurs, avec $E \subset B(H),F
\subset B(K)$ (o\`u $H,K$ sont des Hilbert) on notera $E \otimes_{\min}F$ le
produit tensoriel compl\'et\'e pour la norme induite par l'espace $B(H \otimes K)$
des op\'erateurs born\'es sur le produit tensoriel hilbertien $H \otimes K$. Nous
renvoyons \`a [2] ou [6] pour la d\'efinition du dual $E^*$ d'un espace d'op\'erateur
$E$, ainsi que pour celle de l'espace transpos\'e $E^{op}$. Rappelons simplement
qu'un espace d'op\'erateurs est un espace de Banach $E$ muni d'un plongement
isom\'etrique dans l'espace $B(H)$ des op\'erateurs born\'es sur un Hilbert $H$ et que 
  les
morphismes de la cat\'egorie correspondante sont les op\'erateurs ``compl\'etement
born\'es" (en abr\'eg\'e c.b.) au sens e.g. de [2,5,6,7]. Une application $u : E \to
F$ est compl\`etement born\'ee (resp. compl\`etement isom\'etrique) si l'application
$\tilde u = I_{B(\ell_2)}\otimes u$ s'\'etend en un op\'erateur born\'e (resp.
isom\'etrique) de $B(\ell_2) \otimes_{\min} E$ dans $B(\ell_2) \otimes_{\min} F$
et on n $\N {u}_{cb} = \N {\tilde u}$. On notera
$cb(E,F)$ l'espace des applications c.b. de $E$ dans $F$
muni de cette norme.

Si $E$ et $F$ sont compl\`etement isom\'etriques on note $E \approx F$. On peut
d\'efinir le complexe conjugu\'e $\adh {E}$ d'un espace d'op\'erateurs $E \subset
B(H)$ comme \'etant le m\^eme espace de Banach muni de la multiplication complexe
conjugu\'ee $\lambda.x = \adh {\lambda}x$. On le munit de la structure d'espace
d'op\'erateurs induite par le plongement lin\'eaire $\adh {E} \subset \adh {B(H)}=
B(\adh {H})$. 

\proclaim THEOREME 1. Soit $I$ un ensemble arbitraire. Il existe un espace de
Hilbert ${\cal H}$ et un sous-espace $E$ de $B({\cal H})$ qui est isom\'etrique \`a
$\ell_2(I)$ et tel que l'identification naturelle entre $E$ et $\adh {E^*}$
soit une isom\'etrie compl\`ete. De plus, cet espace est unique \`a une isom\'etrie
compl\`ete pr\`es.

Nous noterons cet espace $OH(I)$ et si $I = \nat$ on note simplement $OH =
OH(\nat)$. C'est l'analogue de l'espace de Hilbert dans la cat\'egorie des
espaces d'op\'erateurs. L'unicit\'e dans le th\'eor\`eme pr\'ec\'edent signifie que tout
espace d'op\'erateur $E$ v\'erifiant les propri\'et\'es du th\'eor\`eme 1, est compl\`etement
isom\'etrique \`a l'espace $OH(I)$. On peut d\'ecrire la structure de $OH(I)$ de la
mani\`ere suivante. Soit $(a_i)_{i \in I}$ une famille \`a support fini d'\'el\'ements
de $B(\ell_2)$ et soit $(T_i)_{i \in I}$ une base orthonormale de $OH(I)$. On a
alors 
$$\N {\sum T_i \otimes a_i}_{OH(I) \otimes_{\min}B(\ell_2)} = \N {\sum a_i
\otimes \adh {a_i}}_{B(\ell_2) \otimes_{\min}B(\adh {\ell_2})}^{1/2}.$$
Il en r\'esulte que l'espace $OH(I)$ est aussi compl\`etement isom\'etrique \`a son
oppos\'e $OH(I)^{opp}$. De plus, tout op\'erateur born\'e $u$ sur $OH(I)$ est
automatiquement c.b. et on a $\N {u}_{cb}= \N {u}$. Soit $E \otimes_h F$ le
produit tensoriel de Haagerup de deux espaces d'op\'erateurs introduit par Effros
et Kishimoto [5] et \'etudi\'e dans [2,6]. On peut montrer que l'on a une
identification compl\`etement isom\'etrique $OH(I) \otimes_h OH(J) = OH(I \times 
J)$.

Rappelons qu'en th\'eorie de l'interpolation [1] on appelle ``compatible" un
couple d'espaces de Banach muni d'un espace vectoriel topologique ${\cal X}$ et
d'injections continues $E_0 \hookrightarrow {\cal X}$ et $E_1 \hookrightarrow
{\cal X}$. Cette structure permet de d\'efinir les espaces $E_0 \cap E_1$ et $E_0
+ E_1$ muni des normes suivantes
$$\forall~x \in E_0 \cap E_1~~~~~~~\N {x}_{E_0 \cap E_1} = \max \{\N
{x}_{E_0},\N {x}_{E_1}\}$$
et
$$\forall~x \in E_0 + E_1~~~~~~~\N {x}_{E_0 + E_1} = \inf \{\N {x_0}_{E_0} +
\N {x_1}_{E_1}|x = x_0 + x_1\}.$$
Signalons que les espaces $R\cap C$ et $R+C$ sont \'etudi\'es dans [8].

\noindent On notera $E_\theta = (E_0,E_1)_\theta$ et $E^\theta =
(E_0,E_1)^\theta$ les espaces d'interpolation complexe introduits par Calder\'on
et Lions (\cf [1]). La caract\'erisation des espaces d'op\'erateurs due \`a Ruan [11]
permet de munir le dual ou bien le quotient d'un espace d'op\'erateurs d'une
structure naturelle d'espace d'op\'erateurs. 
 On note 
$E_0 \oplus_\infty E_1$ (resp. $E_0 \oplus_1 E_1$) la
somme directe munie de la   norme $\|(x,y)\|=
\max(\|x\|,\|y\|)$ (resp. $\|x\|+\|y\|)$. On munit ces
espaces d'une structure d'espace d'op\'erateurs en posant 
  $M_n(E_0 \oplus_\infty E_1)=M_n(E_0
)\oplus_\infty M_n(E_1)$
(resp. $M_n(E_0 \oplus_1 E_1)$ est muni de la norme
induite par la norme de l'espace 
$cb(E_0^*\oplus_{\infty} E_1^*,M_n)$.

 On peut aussi munir les espaces $E_0 \cap E_1$
et $E_0 + E_1$ d'une structure naturelle d'espace
d'op\'erateurs. On identifie  $E_0 \cap E_1$ au sous-espace
$\Delta=\{(x,x)\}$ dans $E_0 \oplus_\infty E_1$ et on munit
$E_0 + E_1$ de la structure d'espace
d'op\'erateurs du quotient
 $(E_0 \oplus_1 E_1)/\{(x,-x)\}$.

On a plus g\'en\'eralement le r\'esultat
suivant.

\proclaim THEOREME 2. Soit $(E_0,E_1)$ un couple compatible d'espaces
d'op\'erateurs. On peut munir  $E_\theta$ et $E^\theta$
d'une structure d'espace d'op\'erateurs en posant 
$$\eqalign{ M_n(E_\theta) & =
(M_n(E_0),M_n(E_1))_\theta~~\cr {\rm (resp.} \quad 
M_n(E^\theta) & = (M_n(E_0),M_n(E_1))^\theta){.}}$$ Si on
suppose $E_0 \cap E_1$ dense \`a la fois dans $E_0$ et dans
$E_1$ on a les identifications compl\`etement
isom\'etriques suivantes : $$(E_0 \cap E_1)^* = E_0^* +
E_1^*,\quad  (E_0 + E_1)^* = E_0^* \cap E_1^*~\quad {\rm
et}~\quad(E_0,E_1)_\theta^* = (E_0^*,E_1^*)^\theta.$$

\proclaim THEOREME 3. Soit $E$ un espace d'op\'erateurs. Soit $v : E \to
\ell_2(I)$ une injection continue \`a image dense de sorte que $J = \adh {v^*} v
: E \to \adh {E^*}$ est une injection continue permettant de consid\'erer
$(E,\adh {E^*})$ comme un couple compatible. On a alors un isomorphisme
compl\`etement isom\'etrique
$$(E,\adh {E^*})_{1/2} \approx OH(I).$$
En particulier $(R,C)_{1/2} \approx OH$.

\proclaim COROLLAIRE 4. Soit $\ell_2^n$ l'espace hilbertien complexe  de
dimension $n$. Consid\'erons les espaces d'op\'erateurs  $R_n=B(\ell_2^n,\comp)$ et
$C_n=B(\comp,\ell_2^n)$ identifi\'es \`a $\ell_2^n$ de la faon \'evidente et notons
$OH_n$ la version $n$-dimensionnelle de $OH$. Soit $H$ un Hilbert. On a alors
compl\`etement isom\'etriquement $$(R_n\otimes_{\min} B(H),C_n\otimes_{\min}
B(H))_{1/2} \approx OH_n\otimes_{\min} B(H).$$

 En utilisant des r\'esultats de Connes [4],
g\'en\'eralis\'es par Haagerup [7] on obtient 

\proclaim COROLLAIRE 5. Une alg\`ebre de von Neumann $M \subset B(H)$ est
injective si (et seulement si) il existe une projection compl\`etement born\'ee
$P$ de $B(H)$ sur $M$.

En fait, il suffit qu'il existe une projection telle que l'on ait pour une
constante $C$ pour toute suite finie dans $B(H)$
$$\N {\sum P(x_i)^* P(x_i)} \leq C^2 \N {\sum x_i^* x_i}~{\rm et}~\N {\sum
P(x_i) P(x_i)^*}\leq C^2 \N {\sum x_i x_i^*}.$$
On peut montrer qu'il suffit m\^eme que l'inclusion naturelle de $M \otimes_{\min}
(R + C)$ dans $M \otimes_{\min} R + M \otimes_{\min} C$ soit born\'ee.
 Soit $E,F$ des espaces d'op\'erateurs. Soit $u \in E \otimes F$. On d\'efinit
$$oh(u) = \inf \{\N {\sum x_i \otimes \adh {x_i}}_{E \otimes_{\min} \adh
{E}}\N {\sum y_i \otimes \adh {y_i}}_{F \otimes_{\min} \adh {F}}\}$$
o\`u l'infimum porte sur toutes les repr\'esentations de $u$ de la forme $u =
\sum_1^n x_i \otimes y_i,\quad  x_i \in E,\quad y_i \in F$.
On v\'erifie ais\'ement que $oh$ est une norme sur $E \otimes F$ et l'on note $E
\otimes_{oh}F$ le compl\'et\'e de $E \otimes F$ pour cette norme. Soit $u : E^*
\to F$ l'op\'erateur associ\'e \`a $u$. On peut d\'ecrire de
faon \'equivalente $oh(u) = \inf \{ \N {A}_{cb} \N {B}_{cb}\}$ o\`u l'infimum
porte sur toutes les factorisations possibles de $u$ de la forme $u = AB$ avec
$B : E^* \to OH$ et $A : OH \to F$ de rangs finis.
\proclaim THEOREME 6. Soit $M,N$ deux alg\`ebres de von Neumann injectives. On a
alors une identification isom\'etrique
$$(M \otimes_h N, M^{op} \otimes_h N^{op})_{1/2} = M \otimes_{oh}N.$$

\underbar {Remarques}. Plus g\'en\'eralement, si $E$ est un espace d'op\'erateurs
arbitraire, on peut consid\'erer pour $0 < \theta < 1$ l'espace d'op\'erateurs
$E(\theta) = (E,E^{op})_\theta$. (Noter que $E(0) = E$ et $E(1) = E^{op})$. Si
$F$ est un autre espace d'op\'erateurs le th\'eor\`eme 2 nous donne 
$$(E \otimes_h F, E^{op} \otimes_h F^{op})_\theta \approx E(\theta) \otimes_h
F(\theta)$$
compl\`etement isom\'etriquement. Mais il faut noter que si $E$ est un sous-espace
ferm\'e de $B(H)$ en g\'en\'eral $E(\theta)$ ne s'identifie pas \`a un sous-espace
ferm\'e de $B(H)(\theta)$.

Dans le corollaire suivant on identifiera $M_n$ et $\adh {M_n^*}$ par
l'isomorphisme lin\'eaire $i : M_n \to \adh {M_n^*}$ d\'efini par $i(x)(y) = tr(x
y^*)$. Rappelons que l'on note $cb(E,F)$ l'espace des applications c.b. de $E$
dans $F$ muni de la norme c.b.. Nous pouvons donc identifier les \'el\'ements de
$cb(M_n,M_n)$ avec ceux de $cb(\adh {M_n^*},\adh {M_n^*})$. On a alors 

\proclaim COROLLAIRE 7. L'application lin\'eaire $\Phi : M_n \otimes M_n \to
CB(M_n,M_n)$ d\'efinie par $\Phi(x \otimes y) (a) = x a y$ \'etablit un
isomorphisme isom\'etrique entre $M_n \otimes_{oh} M_n$ et $(cb(M_n,M_n), cb(\adh
{M_n^*}, \adh {M_n^*}))_{1/2}$.

Soit $M \subset B(H)$ et $N \subset B(K)$ des alg\`ebres de von Neumann injectives.
On peut g\'en\'eraliser le corollaire pr\'ec\'edent dans le cadre des espaces
consid\'er\'es par Effros et Kishimoto [5] et Blecher-Smith [3] des applications
compl\`etement born\'es de $B(H,K)$ dans $B(H,K)$ qui sont $(M',N')$ modulaires. 

Dans une prochaine note nous \'etudions les op\'erateurs qui se factorisent
\`a travers l'espace $OH$ par
 des applications compl\`etement born\'ees.

\magnification\magstep1
 \baselineskip = 18pt
 
 \overfullrule = 0pt
  
 \vfill\eject
 \centerline{\bf References}
 
 \item{[1]} J. Bergh and J. L\"ofstr\"om. Interpolation spaces. An
 introduction. Springer Verlag, New York. 1976.
 
 \item{[2]} D. Blecher and V. Paulsen. Tensor products of operator spaces.
 J. Funct. Anal. 99 (1991) 262-292.
 
  \item{[3]} D. Blecher and R. Smith. The dual of the Haagerup tensor
 product. (Preprint) 1990.

 \item{[4]} A.Connes. Classification of injective factors.
Ann. of Math. 104 (1976) 73-116.
  
 \item{[5]} E. Effros and A. Kishimoto. Module maps and
 Hochschild-Johnson cohomology. Indiana Univ. Math. J. 36 (1987) 257-276.
 
 \item{[6]} E. Effros and Z.J. Ruan. A new approach to operator spaces.
 Canadian Math. Bull.34(1991) 329-337.

 \item{[7]} U. Haagerup. Injectivity and decomposition of completely
 bounded maps in ``Operator algebras and their connection with Topology and
 Ergodic Theory''. Springer Lecture Notes in Math. 1132 (1985) 91-116.
 
\item{[8]} U. Haagerup and G. Pisier. Bounded linear
operators between
 $C^*$-algebras. To appear.

 \item{[9]} O. Kouba. Interpolation of injective or projective tensor
 products of Banach spaces. J. Funct. Anal. 96 (1991) 38-61.
 
 \item{[10]} G.Pisier. The operator Hilbert space and complex interpolation.
(Preprint) A para\^\i tre.

\item{[11]} Z.J. Ruan. Subspaces of $C^*$-algebras. J. Funct. Anal. 76
 (1988) 217-230.
 \medskip

Equipe d'Analyse

Universit\'e Paris 6

Bo\^\i te 186, Place Jussieu

75252 Paris Cedex 05, FRANCE

et 

Mathematics Department

Texas A. and M. University

College Station, Texas 77843, USA
 
 \end

\end